\documentclass[10pt]{article}
\usepackage{amsmath}
\usepackage{amssymb}
\usepackage{latexsym}

\usepackage[all]{xy}

\newtheorem{thm}{Theorem}[section]
\newtheorem{prop}[thm]{Proposition} \newtheorem{lemma}[thm]{Lemma}
 \newtheorem{dfn}[thm]{Definition}
 \newtheorem{rmk}[thm]{Remark}
\newtheorem{ex}[thm]{Example} \newtheorem{question}[thm]{Question}

\newcommand {\pf}{\noindent{\bf Proof.}\ }
\newcommand{\complex}{{\mathbb C}}

\newcommand{\reals}{{\mathbb R}}

\newcommand{\integers}{{\mathbb Z}}

\newcommand{\Hom}{{\rm Hom}}
\newcommand{\End}{{\rm End}}

\newcommand{\aop}{A^{\mathrm{op}}}

\newcommand{\qed}{\begin{flushright} $\Box$\ \ \ \ \ \end{flushright}}

\newcommand{\arrows}{\,\lower1pt\hbox{$\longrightarrow$}\hskip-.24in\raise2pt
          \hbox{$\longrightarrow$}\,}

\newcommand{\bDelta}{\mbox {\boldmath $\Delta$}}
\newcommand{\bPhi}{\mbox {\boldmath $\Phi$}}
\newcommand{\bS}{\mbox {\boldmath $S$}}
\newcommand{\boldf}{\mbox {\boldmath $f$}}
\newcommand{\boldg}{\mbox {\boldmath $g$}}

\newcommand{\bepsilon}{\mbox{\boldmath $\epsilon$}}

\newcommand{\Alg}{\rm \textsf{Alg}}

\begin{document}

\title{{\bf Hopfish algebras}}
\author
{Xiang Tang,
Alan Weinstein\thanks{Research partially supported by NSF Grant
DMS-0204100.\newline
\mbox{~~~~}MSC2000 Subject Classification Number: 16W30(Primary),
81R50(Secondary) .
\newline \mbox{~~~~}Keywords: Hopf algebra, hopfish algebra, groupoid,
bimodule, Morita equivalence, hypergroupoid.}
, and Chenchang Zhu
}
\maketitle
\begin{abstract}
We introduce a notion of ``hopfish algebra'' structure on an
associative algebra, allowing the structure morphisms (coproduct,
counit, antipode) to be bimodules rather than algebra
homomorphisms. We prove that quasi-Hopf algebras are examples of
hopfish algebras. We find that a hopfish structure on the algebra
of functions on a finite set $G$ is
 closely related to a
``hypergroupoid" structure on $G$. The Morita theory of  hopfish algebras is also
discussed.
\end{abstract}
\section{Introduction}
\label{sec-intro}

When the multiplication on a (discrete, topological, smooth, algebraic...)
group $G$ is encoded in an appropriate algebra $A=A(G)$ of
functions on $G$ with values in a commutative ring $k$, it
becomes a coproduct, i.e. an algebra homomorphism
$\Delta:A\to A\otimes_kA.$ The inclusion of
the unit and the inversion map are also encoded as
homomorphisms: the counit $\epsilon:A\to k$ and the antipode
$S:A\to A$. The group properties (associativity, unit,
inverse) become statements about these homomorphisms which
constitute the axioms for a (commutative) {\bf Hopf algebra}; any
noncommutativity of the underlying group appears as
non\underline{co}commutativity of the coproduct.

In noncommutative geometry, a noncommutative algebra
$A$ is thought of as the functions on a ``noncommutative
space'' or ``quantum space'' $X$.  If $X$ is to be a ``quantum
group'', the algebra $A$ should have the additional structure
of a Hopf algebra.  We note that, for noncommutative Hopf algebras, the antipode has to be an
antihomomorphism rather than a homomorphism of algebras.  For this
reason, a Hopf algebra is not quite a group in the category of
algebras; this anomaly will come back to haunt us later.

One type of quantum space is a quantum torus, whose function
algebra is the crossed product algebra $A_\alpha$ associated
to an action of $\integers$ on the circle $S^1=\reals/\integers$
 generated
by an irrational rotation $r_\alpha$.  This irrational rotation
algebra is generally
taken as a surrogate for the algebra of continuous functions on
the ``bad quotient space'' $S^1/\alpha\integers$ because, for nice
quotients, the crossed product algebra is Morita equivalent to the
algebra of functions on the quotient.   Since
$S^1/\alpha\integers$ is a group, one might expect $A_\alpha$
to have a Hopf algebra structure, but this is not so.  In
particular, there can be no counit, since there are no algebra
homomorphisms $A_\alpha \to \complex$.  In geometric language,
``the quantum torus has no points''.

Additionally, in noncommutative geometry, Morita equivalent algebras are often thought of as representing the ``same space'', but the notion of Hopf algebra, and even that of biunital bialgebra, is far from Morita invariant.

In this paper, we propose a new algebraic approach to ``group structure''
 based on the idea that the
appropriate morphisms between algebras are bimodules (perhaps with
extra structure, or satisfying extra conditions) rather than
algebra homomorphisms.  Our immediate inspiration to use
bimodules was the work of Tseng and Zhu \cite{tz}, in which leaf
spaces of foliations are treated as differential stacks for the
purpose of putting group(oid)-like structures on them.  This means
that the structure morphisms of the groupoids are themselves
``bibundles''  \cite{mr:stability} (with respect to foliation groupoids, which play in
this geometric story the role of the crossed product algebras
above) rather than ordinary mappings of leaf spaces.  We were also
motivated by previous uses of bimodules as generalized morphisms of
algebras, $C^*$-algebras, groupoids, and Poisson manifolds, a point of view which has been extensively
developed by Landsman and others (see, for instance, \cite{bu-we:poisson}, \cite{la:bicategories}, and
\cite{la:operator}).

We call our new objects {\bf hopfish algebras}, the
suffix ``oid'' and prefixes like ``quasi'' and ``pseudo'' having
already been appropriated for other uses.  Also, our term retains
a hint of the Poisson geometry which inspired some of our work.

\noindent {\bf Outline of the paper.}  We begin with a discussion of
the category in which objects are algebras and morphisms are
bimodules, emphasizing the functor, which we call modulation, from the
usual category to this one.  We then look at the analogues of
semigroups and groups in this category, which we call sesquialgebras
and hopfish algebras.  What turns out to be especially delicate is the
definition of the antipode.
We next show that Hopf algebras, and the more
general quasi-Hopf algebras, become hopfish algebras upon modulation.
In the following section, we study the Morita invariance of the
hopfish property, showing that a sufficient condition for this to hold
is that a Morita equivalence bimodule is compatible with the antipode
of a hopfish algebra.  Finally, we study hopfish structures on finite
dimensional commutative algebras.  We show that these correspond to
``multiple-valued groupoid structures," and we give examples of
hopfish algebras which do not correspond under Morita equivalence to
Hopf algebras.  

\noindent {\bf Outlook.}  In the present paper, we restrict ourselves
to the purely 
algebraic situation; in particular,  our tensor products do not involve
any completion.  We do not require
finite dimensionality of our algebras, although some of our examples
do have this property.  We hope to develop a theory of hopfish
$C^*$-algebras in the future, with a treatment of irrational rotation
algebras as a first goal.  Even without this theory, two of the
authors, along with Blohmann \cite{bl-ta-we:hopfish}, have succeeded in 
constructing a sesquiunital sesquialgebra structure on the
``polynomial part'' of the irrational rotation algebras.  These
algebras are not quite hopfish, since the candidate antiautomorphism
satisfies only a weakened version of our antipode axiom.  (We hope
that this will be remedied when we go on to the $C^*$-algebras.)
Nevertheless, our structure is sufficient to induce an
interesting monoid structure on isomorphism classes of
modules.  

Finally, we remark that all of our examples of hopfish examples are
either weak Hopf algebras or Morita equivalent to quasi-Hopf algebras.
It would be interesting to find completely new examples.  The
irrational rotation algebras are probably not of either of these
special types, but, as we have already noted, they are not quite hopfish.

\noindent {\bf Acknowledgements.} This work began in July, 2004, when
Tang and Weinstein were participants in the trimester on K-Theory and
Noncommutative Geometry at the Centre Emile Borel.  We had further
opportunities for collaboration, and received further stimulation,
during a number of other conferences and short-term visits. We would
like to thank Max Karoubi, Yvette Kosmann-Schwarzbach, Ryszard Nest,
Tudor Ratiu, Pierre Schapira, and Michel van den Bergh for their
invitations and hospitality.  For useful comments in conversation and
correspondence, we thank Paul Baum, Christian Blohmann,
Micho Durdevic, Piotr Hajac,
Yvette Kosmann-Schwarzbach, Olivier Mathieu, Ryszard Nest, Radu
Popescu, Jean Renault, Earl Taft, and Boris Tsygan. We would
especially like to thank Pavel Etingof for helping us to overcome a
stumbling block in the characterization of the antipode, as well as
Noah Snyder for providing the example in the last section.

\section{The modulation functor}

Fixing a commutative ring $k$  as our ring of
scalars, we will work mostly in a category \Alg~ whose objects are unital
$k$-algebras.  The
morphism space $\Hom(A,B)$ is taken to be the set of
isomorphism classes of biunital $(A,B)$-bimodules. We will
almost always consider these morphisms as going from right to
left, i.e.\- from $B$ to $A$ (or, better, ``to $A$ from $B$'').
The composition $XY \in \Hom(A,C)$ of $X \in \Hom(A,B)$ and $Y
\in \Hom(B,C)$ is defined (on representative bimodules) as
$X\otimes_B Y$, with the residual actions of $A$ and $C$ providing
the bimodule structure.

We will frequently fail to distinguish between morphisms in \Alg~ and their representative bimodules, as long as we can do so without causing confusion.
It is also possible to work in the more refined 2-category whose
morphisms are bimodules and whose 2-morphisms are bimodule
isomorphisms, but we leave this for the future.

We will denote by $\Alg_0$ the ``usual'' category whose objects
are again unital $k$ algebras but whose morphisms are unital
homomorphisms.  Thus, $\Hom_0(A,B)$ will denote the
homomorphisms to $A$ from $B$. There is an important functor from
$\Alg_0$ to $\Alg$ which we will call
{\bf modulation}.\footnote{We are indebted to
Yvette Kosmann-Schwarzbach for
 suggesting this apt name for a functor which is ubiquitous in the
 literature on Morita equivalence, but which does not seem to have
 acquired a standard designation.}  The modulation of $f\in \Hom_0(A,B)$ is the isomorphism class of
$A_f$, which is the $k$-module $A$ with the $(A,B)$-bimodule structure
\begin{equation}
\label{eq-modulation} a\cdot x\cdot b=axf(b).
\end{equation}
We will often denote the modulation of a morphism by the same
symbol, but in bold face, e.g.  $\boldf \in \Hom(A,B)$.
The modulation functor is not necessarily faithful, as the next lemma shows.

\begin{lemma}
\label{lem-conjugate} For  $f,g \in \Hom_0(A, B)$, their
modulations $\boldf$ and $\boldg$ are equal (i.e. the bimodules $A_f$ and $A_g$
are isomorphic) if and only if $f=\phi g\phi^{-1}$  for some
invertible $\phi\in A$.
\end{lemma}
\pf  If $f=\phi g \phi^{-1}$, $\boldf$ and $\boldg$ are
both represented by $A$, with the same left $A$-module structures.
 To correct for the difference between the
right actions of $B$, we introduce the bijective map $\bPhi: A_f\to
A_g$ defined by $x\mapsto x\phi$, which is a bimodule isomorphism  because
$$\bPhi(axf(b))=axf(b)\phi=ax\phi\phi^{-1}f(b)\phi=ax\phi
g(b) = a\bPhi(x)g(b).$$

For the converse,  given a bimodule isomorphism
$\bPhi:A_f \to A_g,$  we define $\phi$ to be $\bPhi(1_A)$.
 By setting $x=1_A$ in the bimodule
morphism identities $\bPhi(ax)=a\bPhi(x)$ and
$\bPhi(xf(b))=\bPhi(x)g(b)$, we find first that  $\bPhi(a)=a\phi$, so that $\phi$ is invertible because $\bPhi$ is, and then that
$f(b) \phi = \phi g(b)$, or $f=\phi g \phi^{-1}.$ \qed

\begin{lemma}\label{lem-freerank1}
A morphism $X\in \Hom (A, B)$ is the modulation of $f \in
\Hom_0(A, B)$ if and only if it is isomorphic to $A$ as a left $A$ module.
\end{lemma}
\pf If $X$ represents $f$, then clearly $X$ is isomorphic to $A$
as a left $A$ module. For the converse, if $X=A$ as a left $A$
module then $X$ is isomorphic to $A_f$ where $f(b) = 1_A\cdot b$.
\qed

An invertible morphism in $\Hom(A,B)$ is called a {\bf Morita
 equivalence} between $A$ and $B$, and the group of Morita
 self-equivalences of $A$ is called its {\bf Picard group}.
 The modulation functor clearly
 takes algebra isomorphisms to Morita equivalences.  In fact, we
 have:

\begin{lemma}
\label{lem-invertible} The modulation of $f \in \Hom_0(A,B)$ is
invertible if and only if $f$ is invertible.
\end{lemma}

\pf It is a standard fact about Morita equivalence that, if $X\in
\Hom(A,B)$ is invertible, the natural homomorphisms from $A$ and $B$ to
the $B$- and $A$-endomorphisms of $X$ are isomorphisms.   When
$X=A_f,$ the map which takes $b\in B$ to the operator of right
multiplication by $f(b)$ is injective if and only if $f$ is
injective.  On the other hand, all of the left $A$-module
endomorphisms of $A$ are the right multiplications, so they
are all realized by the action of $B$ if and only if $f$ is
surjective. \qed

\begin{rmk} \label{rmk-nonunital}
{\em It is also possible to modulate a nonunital $f$.  In this
case, the underlying $k$-module should be taken to be the left ideal
$I$ in $A$ generated by $f(1_B)$, so that  the bimodule structure
\eqref{eq-modulation} is still biunital. The three lemmas above
change to the following statements, whose proofs are similar,
so  we only  sketch them. }
\end{rmk}

\noindent
{\bf Lemma \ref{lem-conjugate}$^\prime$}  If $f$ and $g$ are algebra homomorphisms
$A\leftarrow B$ not necessarily unital, then their modulations
$\boldf$ and $\boldg$ are equal if and only if there are $\phi\in A\cdot
f(1_B)$ and $\psi \in A \cdot g(1_B)$ such that $\phi \psi =
g(1_B)$, $\psi \phi = f(1_B)$,
$g=\phi f \psi$, satisfying the two additional conditions that
 $x \phi \psi=0$ implies $x \phi=0$ and $x\psi \phi=0$ implies $x \psi =0$.

\noindent
Sketch of the proof: Given an isomorphism $\Phi$ to
$\boldf$ from $\boldg$,  let
$\phi =\Phi(g(1_B))$ and $\psi = \Phi^{-1}(f(1_B))$. Then $\Phi(x
g(1_B)) =x\phi$ and $\Phi^{-1}(x f(1_B))=x\psi$. All this gives us
the desired equations and properties. For the converse, the
morphism $\Phi(x g(1_B)) := x\phi$ is an isomorphism from $A \cdot
g(1_B)$ to $A\cdot f(1_B)$ with inverse $\Phi^{-1}( x f(1_B)):=
x\psi$. The two additional conditions make $\Phi$ and $\Phi^{-1}$
well defined.  \qed

\noindent {\bf Lemma \ref{lem-freerank1}$^\prime$}  A morphism
$X\in \Hom (A, B)$ is the modulation of a (not necessarily unital) map
$f:A\leftarrow B$ if and only if it is represented by a principal left
ideal in $A$.

\noindent Sketch of the proof: If $X$ is the modulation of $f$,
then $X=A\cdot f(1_B)$. For the converse, if $X$ is isomorphic to
a left $A$ ideal $A\cdot c$, then $X$ is the modulation of $f:
b\mapsto c\cdot b$, where $b\in B$ and ``$\cdot$'' is the right
action of $B$ on $X=A\cdot c$.  \qed

\noindent {\bf Lemma \ref{lem-invertible}$^\prime$}   When
$f(1_B)$ is in the center of $A$, the modulation of a morphism $f:
A\leftarrow B$ (not necessarily unital) is invertible if and only if $f$ is
an isomorphism from $B$ to $A\cdot f(1_B)$ and $f(1_B)$ is not a
zero divisor .

\noindent Sketch of the proof: One applies the same argument. If
$\boldf$ is invertible, notice that $A \to \End_B(X)$ by $a\mapsto
a\cdot$ is an isomorphism, therefore $af(1_B)\neq a' f(1_B)$ if
$a\neq a'$. This implies that $f(1_B)$ is not a zero divisor. As
before $f$ has to be injective. For any $a\in A$, right
multiplication by $a$ is in $\End_A(X)$, therefore there is $b\in
B$ such that $f(1_B) a= f(b)$. It is not hard to prove the converse.\qed

Finally, we recall that every $(A,B)$ bimodule gives
rise (via tensor product over $B$) to a $k$-linear functor from
the category of left $B$-modules to that of left $A$-modules, that
isomorphisms between bimodules produce naturally equivalent
functors, and that invertible elements of $\Hom(A,B)$ correspond
to homotopy classes of equivalences of categories.  (The Eilenberg-Watts theorem characterizes the functors arising from bimodules as those which commute with finite limits and colimits.)

\subsection{Sesquialgebras}

To make the notion of biunital bialgebra Morita invariant,
we introduce the following definition.  For simplicity of notation, we
omit the subscript $_k$ on tensor products over $k$, and the unadorned asterisk $^*$ will denote the $k$-dual.

\begin{dfn}
\label{dfn-sesqui} A {\bf sesquiunital sesquialgebra} over a
commutative ring $k$ is a unital $k$-algebra $A$ equipped with
an $(A\otimes A, A)$-bimodule $\bDelta$ (the {\bf
 coproduct}) and a $(k,A)$-module (i.e.
 a right $A$ module) $\bepsilon$
(the {\bf counit}), satisfying the following properties.

\begin{enumerate}
\item (coassociativity)  The
 $(A\otimes A\otimes A,A)$-bimodules
$(A\otimes \bDelta)\otimes_{A\otimes  A}\bDelta$ and
$(\bDelta\otimes A)\otimes_{A\otimes  A}\bDelta$ are
isomorphic.
\item (counit)  The $(k\otimes A,A)=
(A\otimes  k,A)=(A,A)$-bimodules \\
 $(\bepsilon\otimes  A)\otimes_{A\otimes  A}\bDelta$ and
$(A\otimes \bepsilon)\otimes_{A\otimes A}\bDelta$
are both
 isomorphic to $A$.
\end{enumerate}
\end{dfn}

For example, if $(A, \Delta, \epsilon)$ is a biunital bialgebra,
then its modulation  $(A, \bDelta, \bepsilon)$ is a sesquiunital
sesquialgebra.  If we have a Morita equivalence $X$ between $A$
and another algebra $B$, we can use composition with $X$ and
$X\otimes X$ to put a biunital sesquialgebra structure on $B$.
See Section \ref{sec-morita} below for more details.

\section{The antipode and hopfish algebras}
Our definition of sesquiunital sesquialgebra expresses (with
arrows reversed) the usual axioms of a monoid (semigroup with
identity) in the category \Alg. A monoid is a group when all its
elements have inverses, so it is natural to look for a
sesquialgebraic analogue of the inverse.  In a Hopf algebra, the
antipode, which encodes inversion, is an algebra {\em
anti}homomorphism $S:A\to A$.  The properties of inversion
($gg^{-1}=e=g^{-1}g$ for every group element) are then expressed
as commutativity of two diagrams, or equality of compositions

\begin{equation}
\label{eq:antipode} 1\circ\epsilon = \mu\circ\beta\circ\Delta,
\end{equation}
where $1:k\to A$ is inclusion of the scalars,
$\mu:A\otimes A\to A$ is algebra multiplication, and
$\beta:A\otimes A\to A\otimes A$ is either
$I\otimes S$ or $S \otimes I$ ($I$ being the identity morphism on
$A$).

When $A$ is noncommutative, the maps $\mu$ and $\beta$ are
$k$-linear but not algebra homomorphisms.  One can consider $S$ as
a homomorphism from $A$ to the opposite algebra
$A^{\mathrm{op}}$, or vice versa, but there is no way to
correct $\mu$ in such a manner.  As a result, we see no way to
rewrite \eqref{eq:antipode} in the category \Alg.  Instead, we
take an alternate approach,
 which may also be useful elsewhere in the theory of Hopf algebras.

We keep in mind the example where $A$ is the algebra of $k$-valued functions on a group $G$.

One way to characterize groups among monoids without explicitly postulating the existence of inverses is to consider the subset
$$J=\{(g,h)|gh=e\} \subset G \times G$$ and require that it project
 bijectively to one factor in the product.
To represent $J$ algebraically, even when $A$ is noncommutative,
we borrow an idea from Poisson geometry \cite{lu:mmap}, where
 coisotropic submanifolds become one-sided ideals
when a Poisson manifold is quantized to become a noncommutative
algebra.

We begin, then, with the space $Z'=\Hom_A(\bepsilon,\bDelta)$
of right module homomorphisms. (In the group case, $Z'$ plays the role of
``measures" on $G\times G$ which are supported on  $J$.)
Using the left
$A\otimes A$ module structure on $\bDelta$, we define a
right $A\otimes A$ module structure on $Z'$ by
$(gb)(u)=g(bu)$ for $g$ in $Z'$, $b$ in $A\otimes A$ and
$u$ in $\bDelta$.  Note that  $Z'$ is completely determined by $\bepsilon$ and
$\bDelta$ and is {\em not} an extra piece of data.

For the algebraic model of functions on $J$, we must take a {\bf predual} of $Z'$, i.e. a left $A\otimes A$-module $Z$ whose $k$-dual $Z^*$ is equipped with a right $A\otimes A$-module isomorphism with $Z'$.

\begin{dfn}\label{dfn-preantipode}
A {\bf preantipode}
 for a sesquiunital sesquialgebra $A$ over
$k$ is a left $A\otimes  A$ module $\bS$ together with an
isomorphism of its $k$-dual with the right $A\otimes A$
module $\Hom_A(\bepsilon,\bDelta)$.
\end{dfn}

 Since a left $A$ module is
also a right $A^{\mathrm{op}}$ module, we may consider $\bS$ as an
$(A,\aop)$ bimodule, where $(A,\cdot)$ is from the left $A$ in
$A\otimes A$ and $(\cdot,\aop)$ is from the right one, i.e. as an
$\Alg$~morphism in $\Hom(A,\aop)$.

The following is our way of expressing algebraically that the first projection from $J$ to $G$ is bijective.

\begin{dfn}
\label{dfn-antipode}
Let $A$ be a sesquiunital sesquialgebra.
If a preantipode $\bS$, considered as an $(A, A^{op})$ bimodule, is a
free left $A$ module of rank 1, we call $\bS$ an {\bf antipode}
and say that $A$ along with $\bS$ is a {\bf hopfish
algebra}.
\end{dfn}

By Lemma \ref{lem-freerank1}, $\bS$ is the modulation of an
algebra homomorphism $A\leftarrow A^{op}$.   Thus, the definition
is effectively that there is a homomorphism $S$ to $A$ from
$A^{op}$ such that the full $k$-dual of the modulation of $S$ is
isomorphic to $\Hom_A(\bepsilon, \bDelta)$.

\section{Hopf and quasi-Hopf algebras as hopfish algebras}

As we observed earlier, the modulation of a biunital bialgebra is
a sesquiunital sesquialgebra.  In this section, we will give an
explicit description of a pre-antipode in this case, and we will
show that the modulation of a Hopf algebra is hopfish.  Although
this is a special case of the quasi-Hopf algebras treated in the
next section, we deal separately with the Hopf case because the
proof is much simpler.

Let $(A,\Delta,\epsilon)$ be a biunital bialgebra.  Considering
the modulations $\bepsilon=k$ and $\bDelta=A\otimes A$ as
right $A$ modules respectively, one may identify $Z'$ with the
subspace of  $(A\otimes
A)^*=\Hom_k(k,A\otimes A)$ consisting of those linear
functionals which annihilate the left ideal $W$ generated by
$$\{\epsilon(a)(1\otimes
1)-\Delta(a)|a\in A\},$$ i.e. with the $k$-module dual to
$(A\otimes  A)/W$.  We may therefore take the (cyclic)
left $A\otimes A$ module $\bS_1=(A\otimes A)/W$ as
a preantipode.

We will use the following lemma later.  Its straightforward proof is left to the reader.

\begin{lemma}
\label{lemma-delkereps}
 {\em $W$ is equal to the left ideal generated by
$\Delta(\ker\epsilon).$ }
\end{lemma}

Now suppose that $A$ is equipped with an antipode $S$ making it
into a Hopf algebra.  We will consider $S$ as a homomorphism
$A\leftarrow A^{op}$, with modulation $\bS$.  As a $k$-module,
$\bS$ is $A$; its  $(A,A^{op})$ bimodule structure is $ a\cdot x
\cdot b =  ax S(b). $

If  we can show that the preantipode $\bS_1$ is isomorphic to
$\bS$ as a bimodule, then since $\bS$ is isomorphic to $A$  as a
left $A$-module,  $\bS=\bS_1$ is an antipode, making the
modulation of $A$ into a hopfish algebra.

We define a map $\phi: A\otimes A \to A$ by
\[  a \otimes b \mapsto a S(b),\]
This map is obviously a morphism of $(A, A^{op})$ bimodules
because
\[ \phi(c \cdot  (a\otimes b))=\phi( c a\otimes b)=  c a S(b)=
c\cdot(a S( b)) , \] and
\[ \phi((a\otimes b)\cdot c)=\phi(a\otimes cb))= a S( b) S(c) = (a
S(b))\cdot c.\] Hence this map descends to $\bS_1=(A\otimes A)/W$
because
\[ \phi(\epsilon(a)(1\otimes 1)-\Delta(a))= 1\circ\epsilon(a)-(Id\otimes
S)\circ \Delta (a)=0.\] The induced map from $\bS_1$ to $A$, which
we also denote by $\phi$, is also a morphism of $(A, A^{op})$
bimodules.

Moreover $\phi$ is surjective, since it has a left inverse $a
\mapsto [a\otimes 1]$, where $[~]$ denotes the equivalence class
modulo $W$.  This map is also a right inverse, and $\phi$ is
injective, if and only if the equation
\begin{equation} \label{check}
 1\otimes a -
S(a)\otimes 1 \in W \end{equation} is satisfied for all $a\in A$.
Notice that  $aS(b)\otimes 1 -a\otimes b= (a\otimes
1)(S(b)\otimes1-1\otimes b)$ and W is a left ideal.  Since $id=
m\circ (id\otimes \epsilon) \circ \Delta$, composing with $S$ we
have $\sum S(a_1)\epsilon(a_2)=\sum S(a_1\epsilon(a_2)) =S(a)$.
(Here we use  Sweedler's notation $\Delta(a)=\sum a_1\otimes a_2$
and $\sum \Delta(a_1)\otimes a_2 =\sum a_{1,1} \otimes a_{1,2}
\otimes a_2$, etc.) On the other hand, we have
\[
\begin{split}
 & \sum (S(a_1)\otimes 1) \cdot \Delta (a_2)=\sum (S(a_1)a_{2,
1})\otimes a_{2,2}\\
=&\sum (S(a_{1, 1}) a_{1, 2})\otimes a_{2}= \sum 1\otimes
\epsilon(a_1) a_2=1\otimes a.
\end{split}
\]
We explain the equalities above as follows.   The first equality
just comes from the notation and the multiplication in the tensor
product algebra.  For the second, we consider the map $s:A\otimes
A\otimes A\to A\otimes A$ defined by $s(a\otimes b\otimes
c)=S(a)b\otimes c$. Coassociativity and evaluation of $s$ give
\[
\begin{split}
&\sum s(a_{1}\otimes a_{2,1}\otimes a_{2,2})=\sum s(a_{1,1}\otimes
a_{1,2}\otimes a_2)\\
=&\sum (S(a_1)a_{2, 1})\otimes a_{2,2}=\sum (S(a_{1, 1}) a_{1,
2})\otimes a_{2}.
\end{split}
\]
For the third equality, we have used the property of $S$ that $\mu
\circ (S\otimes id) \circ \Delta= 1\circ \epsilon$. Therefore,
\[1\otimes a- S(a)\otimes 1 =\sum (S(a_1)\otimes 1
)(-\epsilon(a_2)+\Delta(a_2)) \in W.\] So \eqref{check} is proven,
hence $\bS\cong \bS_1$ as $(A, A^{op})$ bimodules.

We have thus proved the following theorem.
\begin{thm}
\label{thm-hopfhopfish} Let $(A,\Delta,\epsilon)$ be a biunital
bialgebra.  Then  $(A\otimes  A)/W$, where $W$ is the left ideal
generated by
 $$\{\epsilon(a)(1\otimes 1)-\Delta(a)|a\in A\},$$ is a preantipode for the modulation of $A$.  If $A$ is a Hopf algebra, with antipode $S$, then $(A\otimes  A)/W$ is isomorphic to the modulation  $\bS$, and $(A,\bDelta,\bepsilon,\bS)$
 is a hopfish algebra.
\end{thm}

\begin{rmk}
{\em The hopfish antipode $\bS$ is also isomorphic to $A^{op}$ as a right
$A^{op}$-module if and only if the Hopf antipode $S$ is
invertible.  This is why we use a ``one sided''
criterion for a preantipode to be an antipode.}
\end{rmk}

We turn now to quasi-Hopf algebras.  Recall that a quasi-bialgebra  $(A, \epsilon, \Delta,
S)$ is nearly a bialgebra, except that the coproduct does not
satisfy associativity exactly; instead, there is an invertible
element $\Phi \in A\otimes A \otimes A$ (the coassociator),
satisfying
\begin{equation}\label{coas}
(id\otimes \Delta)(\Delta(a))=\Phi^{-1}(\Delta\otimes
id)(\Delta(a))\Phi, \quad \forall a \in A,
\end{equation}and further coherence conditions,
\begin{equation}\label{coas-cohe}
\begin{split}
(\Delta\otimes id \otimes id)(\Phi)\cdot(id\otimes id\otimes \Delta)(\Phi) &=(\Phi\otimes 1)\cdot(id\otimes\Delta\otimes id)(\Phi)\cdot(1\otimes \Phi), \\
(\epsilon\otimes id)\circ \Delta= & id = (id\otimes \epsilon)\circ \Delta, \\
 (id\otimes \epsilon \otimes id)(\Phi)& =1.
\end{split}
\end{equation}
Since the modulation functor ``kills'' inner automorphisms (Lemma
\ref{lem-conjugate}), the modulation of a quasi-bialgebra is a
sesquialgebra.

Now $A$ is a quasi-Hopf algebra if there is an anti-homomorphism
$S: A\to A$ and elements $\alpha$, $\beta$ in $A$, such that
\begin{equation}\label{quasi-anti}
\sum S(a_1)\alpha a_2=\epsilon(a)\alpha, \quad \sum a_1\beta
S(a_2)=\epsilon(a)\beta, \quad \forall a\in A,
\end{equation}
where we use Sweedler's notation: $\Delta(a)=\sum a_1\otimes a_2$.
There are also higher coherence conditions for $\alpha$ and
$\beta$, regarding which we refer to \cite{dr:quasi-hopf} for
details.

The following proposition is a slight modification of
Proposition 1.5 of Drinfel'd \cite{dr:quasi-hopf}.   Unlike Drinfel'd, we
do not assume that $S$ is invertible, so we can not obtain the
``right'' part of his proposition, but his ``left'' part can
be proven under weaker hypotheses.

\begin{prop}
\label{prop:quasi-hopf} Let $(A, \Delta, \epsilon, \Phi, S,
\alpha, \beta)$ be a quasi-Hopf algebra, with $\Phi=\sum_i
X_i\otimes Y_i \otimes Z_i$ and $\Phi^{-1}=\sum_j P_j\otimes
Q_j\otimes R_j$. Define $\omega\in A\otimes A$ as $\sum_j
S(P_j)\alpha Q_j\otimes R_j$. Denote by $W$ the left ideal of
$A\otimes A$ generated by $\Delta(\ker \epsilon)$. Then
\begin{enumerate}
\item the $k-$linear mappings $\phi, \psi:A\otimes A\to A\otimes A$,
given by $\phi(a\otimes b)=(a\otimes 1)\omega \Delta(b)$ and
$\psi(a\otimes b)=\sum_i aX_i\beta S(Y_i)S(b_1)\otimes b_2 Z_i$,
are bijective, where we have used Sweedler's notation $\Delta
b=b_1\otimes b_2$.
\item the mapping $a\otimes b\mapsto (id\otimes
\epsilon)(\phi^{-1}(a\otimes b))$ induces a bijection $(A\otimes
A)/W\to A$, and $(id\otimes \epsilon)(\phi^{-1}(a\otimes
b))=a\beta S(b)$;
\end{enumerate}
\end{prop}
\pf First, we prove
\[
\begin{array}{ll}
1)\ \phi\psi=id& 2)\ \psi\phi=id.
\end{array}
\]

In the following, we will only prove  Equation 1) above. One
can prove Equation 2) by the same method, as in
\cite{dr:quasi-hopf}. We have
\[
\begin{split}
&\phi\psi(a\otimes b)\\
=&\sum_i\phi(aX_i\beta S(Y_i)S(b_1)\otimes b_2Z_i)\\
=&\sum_i(aX_i\beta S(Y_i)S(b_1)\otimes 1)\omega \Delta(b_2)\Delta(Z_i)\\
=&\sum_{i}(a\otimes 1)(X_i\beta S(Y_i)\otimes 1)\big((S(b_1)\otimes 1)\omega \Delta b_2\big) \Delta Z_i\\
=&\sum_i (a\otimes 1)(X_i\beta S(Y_i)\otimes 1)(B)\Delta(Z_i),
\end{split}
\]
where $B=(S(b_1)\otimes 1)\omega \Delta b_2$.

We insert the definition of $\omega$ in $B$, and have
\[
\begin{split}
&(S(b_1)\otimes 1)\omega \Delta b_2\\
=&\sum_j \big(S(b_1)S(P_j)\alpha Q_j\otimes R_j\big)\Delta b_2\\
=&\sum_j (m\otimes id)\big( (S\otimes \alpha\cdot \otimes id)\big((P_jb_1\otimes Q_j\otimes R_j)(1\otimes \Delta b_2) \big)\big)\\
=&\sum_j (m\otimes id)\big( (S\otimes \alpha\cdot \otimes id)\big((P_j\otimes Q_j\otimes R_j)(b_1\otimes 1\otimes 1)(1\otimes \Delta b_2)\big) \big)\\
=&(m\otimes id)\big( (S\otimes \alpha\cdot \otimes
id)\big(\Phi^{-1}(1\otimes \Delta)\Delta(b)\big) \big),
\end{split}
\]
where $m:A\otimes A\to A$ is the multiplication on $A$, and
$\alpha\cdot: A\otimes A$ is the left multiplication by $\alpha$.

Using the twisted coassociativity $(id\otimes
\Delta)\Delta=\Phi(\Delta\otimes id)(\Delta) \Phi^{-1}$ we continue the
calculation above to find that $B$ is equal
to
\[
\begin{split}
&(m\otimes id)\big((S\otimes \alpha\cdot \otimes id)\big( (\Delta\otimes id)\Delta(b)\Phi^{-1}\big) \big)\\
=&\sum_j (m\otimes id)\big((S\otimes \alpha\cdot \otimes id)\big( b_{11}P_j\otimes b_{12}Q_j\otimes b_2R_j \big)\big)\\
=&\sum_j (m\otimes id)\big( S(P_j)S(b_{11})\otimes \alpha b_{12}Q_j\otimes b_2R_j\big)\\
=&\sum_j S(P_j)S(b_{11})\alpha b_{12} Q_j\otimes b_2R_j\\
=&\sum_j S(P_j)\alpha \epsilon(b_1)Q_j \otimes b_2R_j\\
=&\sum_j S(P_j)\alpha Q_j\otimes \epsilon(b_1)b_2 R_j\\
=&\sum_j S(P_j)\alpha Q_j \otimes b R_j\\
=&(1\otimes b)\sum_j(S(P_j)\alpha Q_j \otimes R_j),
\end{split}
\]
where in the fourth equality we have used a property of the antipode
$S$, and at the fifth we have used a property of
$\epsilon$.

Substituting the expression above for $B$ in the calculation
of $\phi\psi$, we have
\[
\begin{split}
&\phi\psi(a\otimes b)\\
=&\sum_{i,j}(a\otimes 1)(X_i\beta S(Y_i)\otimes 1)(1\otimes b)(S(P_j)\alpha Q_j\otimes R_j)\Delta (Z_i)\\
=&(a\otimes b)\sum_{i,j}(X_i\beta S(Y_i)\otimes 1)(S(P_j)\alpha
Q_j\otimes R_j)\Delta(Z_i) .
\end{split}
\]

Next, we show that $U=\sum_{i,j}(X_i\beta S(Y_i)\otimes
1)(S(P_j)\alpha Q_j\otimes R_j)\Delta(Z_i)$ is equal to 1. We
introduce the $k$-linear  map $\Psi: A\otimes A\otimes A
\otimes A\to A\otimes A$ by $\Psi(a\otimes b\otimes c\otimes
d)=a\beta S(b) \alpha c\otimes f$, so that $U$ can be written as
\[
\begin{split}
&\sum_{i,j}(X_i\beta S(Y_i)\otimes 1)(S(P_j)\alpha Q_j\otimes R_j)\Delta(Z_i)\\
=&\sum_{i,j}X_i\beta S(Y_i)S(P_j)\alpha Q_jZ_{i1}\otimes R_jZ_{i2}\\
=&\sum_{i,j}\Psi((1\otimes P_j\otimes Q_j\otimes R_j)(X_i\otimes Y_i\otimes Z_{i1}\otimes Z_{i2}))\\
=&\Psi\big((1\otimes \Phi^{-1})(id\otimes id\otimes
\Delta)(\Phi)\big).
\end{split}
\]

Using the coherence condition
\begin{equation}
\label{eq:triple} (id\otimes id\otimes \Delta)(\Phi)(\Delta\otimes
id \otimes id)(\Phi)=(1\otimes \Phi)(id\otimes \Delta \otimes
id)(\Phi)(\Phi\otimes 1),
\end{equation}
we get
\[
\begin{split}
(1\otimes \Phi^{-1})(id\otimes id\otimes \Delta)(\Phi)&=(id\otimes \Delta \otimes id)(\Phi)(\Phi\otimes 1)(\Delta \otimes id\otimes id)(\Phi^{-1})\\
=&\sum_{i,j,k}X_iX_jP_{k1}\otimes Y_{i1}Y_jP_{k2}\otimes
Y_{i2}Z_jQ_k\otimes Z_iR_k.
\end{split}
\]

Hence $\Psi\big((1\otimes \Phi^{-1})(id\otimes id\otimes
\Delta)(\Phi)\big)$ is equal to
\[
\begin{split}
&\sum_{i,j,k}\Psi(X_iX_jP_{k1}\otimes Y_{i1}Y_jP_{k2}\otimes Y_{i2}Z_jQ_k\otimes Z_iR_k)\\
=&\sum_{i,j,k}X_iX_jP_{k_1}\beta S(P_{k2})S(Y_j)S(Y_{i1})\alpha Y_{i2}Z_jQ_k\otimes Z_iR_k\\
=&\sum_{i,j,k}X_iX_j\beta\epsilon(P_k)S(Y_j)\epsilon(Y_i)\alpha Z_jQ_k\otimes Z_iR_k\\
=&\sum_{i,j,k}X_i\epsilon(P_k)(X_j\beta S(Y_j)\alpha Z_j)\epsilon(Y_i)Q_k\otimes Z_iR_k\\
=&\sum_{i,k}X_i\epsilon(P_k)\epsilon(Y_i)Q_k\otimes Z_iR_k.
\end{split}
\]
Here, in the second equality, we have used properties of  the
antipode, i.e. $P_{k1}\beta S(P_{k2})=\beta \epsilon (P_k)$, and
$S(Y_{i1})\alpha Y_{i2}=\alpha \epsilon(Y_i)$. In the last
equality, we have used $\sum_j X_j\beta S(Y_j)\alpha Z_j=1$.

We evaluate $id\otimes \epsilon\otimes id\otimes id$ on both sides
of \eqref{eq:triple}, and since $\epsilon$ is an
homomorphism from $A$ to $k$, we obtain
\begin{equation}
\label{eq:epsilon-id}
\begin{split}
&(id\otimes \epsilon\otimes \Delta)(\Phi)((id \otimes \epsilon)\Delta\otimes id\otimes id)(\Phi)\\
=&(id\otimes((\epsilon\otimes id\otimes
id)(\Phi)))(id\otimes(\epsilon \otimes id)\Delta\otimes
id)(\Phi)((id\otimes \epsilon \otimes id)(\Phi)\otimes id).
\end{split}
\end{equation}

In the definition of a quasi-Hopf algebra, we have assumed that
$id\otimes \epsilon \otimes id(\Phi)=1$. Therefore, $(id\otimes
\epsilon \otimes \Delta)(\Phi)=(id \otimes id\otimes
\Delta)(id\otimes \epsilon \otimes id)(\Phi)=1$. Hence,  by
$(id\otimes \epsilon))\Delta=id\otimes 1$, the left hand side of
(\eqref{eq:epsilon-id}) is equal to
\[
((id\otimes \epsilon)\Delta\otimes id \otimes id)(\Phi)=\sum_i X_i\otimes 1\otimes Y_i\otimes Z_i.
\]
The right hand side of \eqref{eq:epsilon-id} is equal to
\[
(\epsilon\otimes id\otimes id)(\Phi)(\sum\limits_i X_i\otimes 1\otimes Y_i\otimes Z_j).
\]
Therefore, we have
\begin{equation}
\label{eq:eq-8}
\sum_i X_i\otimes 1\otimes Y_i\otimes Z_i=(\epsilon\otimes id\otimes id)(\Phi)(\sum\limits_i X_i\otimes 1\otimes Y_i\otimes Z_j).
\end{equation}

We multiply  both sides of Equation (\ref{eq:eq-8}) by $\sum\limits_i P_j\otimes 1\otimes Q_j\otimes R_j$ and obtain
\[
\epsilon\otimes id\otimes id(\Phi)=1.
\]
So we have $\epsilon \otimes id \otimes id(\Phi^{-1})=\epsilon
\otimes id \otimes id(\Phi^{-1}\Phi)=1$.

Finally, $\sum_{i,j}X_i\epsilon(P_k)\epsilon(Y_i)Q_k\otimes
Z_iR_k$ can be written as
\[
\begin{split}
=&\sum_{i,k}(m\otimes id)(X_i\otimes \epsilon(Y_i)\otimes Z_i )(\epsilon(P_k)\otimes Q_k\otimes R_k)\\
=&(m\otimes id)\big( (id\otimes \epsilon\otimes id)(\Phi)(\epsilon\otimes id\otimes id)(\Phi^{-1})\big)\\
=&1.
\end{split}
\]
In conclusion, we have shown that $\phi\psi(a\otimes b)=a\otimes
b$ and similarly $\psi\phi(a\otimes b)=a\otimes b$. Therefore,
$\phi$ and $\psi$ are invertible.  This  completes the proof of the first
statement of Proposition \ref{prop:quasi-hopf}.

Now calculate $(id\otimes\epsilon)\phi^{-1}(a\otimes b)$. By the
proof above, $\psi$ is the inverse of $\phi$, and
\[
\begin{split}
&(id\otimes \epsilon)\phi^{-1}(a\otimes b)\\
=&(id\otimes \epsilon)(\sum_i aX_i\beta S(Y_i)S(b_1)\otimes b_2Z_i)\\
=&\sum_i aX_i\beta S(Y_i)S(b_1)\epsilon (b_2)\epsilon(Z_i)\\
=&\sum_i aX_i\beta S(Y_i)S(b_1\epsilon (b_2))\epsilon(Z_i)\\
=&\sum_i aX_i\beta S(Y_i)S(b)\epsilon(Z_i).
\end{split}
\]

To show that the last term is is equal to $a\beta S(b)$, we
consider the $k$-linear map $\Upsilon:A\otimes A\otimes A\to A$
defined by $\Upsilon(a_1\otimes a_2\otimes a_3)=a_1\beta
S(a_2)a_3$. Accordingly, we have $\sum_i X_i\beta
S(Y_i)\epsilon(Z_i)=\Upsilon((id \otimes id \otimes
\epsilon)(\Phi))$. By applying $id\otimes id \otimes \epsilon
\otimes id$ to \eqref{eq:epsilon-id}, similarly we have
$(id\otimes id\otimes \epsilon)(\Phi)=1\otimes 1\otimes 1$. So
$\sum_i X_i\beta S(Y_i)\epsilon(Z_i)=\Upsilon(1)=\beta$, and
$\sum_i aX_i\beta S(Y_i)S(b)\epsilon(Z_i)$ is equal to $a\beta
S(b)$.

Therefore if there is an element in $W$, which can be written as
$\Delta(\mu)$, where $\mu$ is in the kernel of $\epsilon$.
$(id\otimes \epsilon)\phi^{-1}(\Delta(\mu))=\mu_1\beta
S(\mu_2)=\epsilon(\mu)\beta=0$. This shows that $W$ is contained
in the kernel of the map $(id\otimes \epsilon)\phi^{-1}:A\otimes
A\to A$. We finally show that $(id\otimes \epsilon)\phi^{-1}$ is a
bijection from $A\otimes A/W$ to $A$. If $\sum_i x_i \otimes y_i
$ is in the kernel of $(id\otimes \epsilon)\phi^{-1}$. We define
$\sum_j a_j\otimes b_j$ equal to $\phi(\sum_i x_i\otimes y_i)$,
and $(id\otimes \epsilon)(\sum_j a_j\otimes b_j)=\sum_j
a_j\epsilon (b_j)=0$. Then $\sum_{i}x_i\otimes y_i$ is equal to
\[
\begin{split}
&\sum_i x_i\otimes y_i\\
=&\sum_j \phi(a_j\otimes b_j)\\
=&\sum_j (a_j\otimes 1)\omega \Delta(b_j)\\
=&\sum_j (a_j\otimes 1)\omega (\Delta(b_j)-\epsilon(b_j))\in W,
\end{split}
\]
where at the third equality, we have used that
$$\sum_{j}(a_j\otimes 1)\omega
\epsilon(b_j)=\sum_j(a_j\epsilon(b_j)\otimes 1)\omega=0.$$ \qed

By the same arguments as in Theorem \ref{thm-hopfhopfish}, we have the following result.
\begin{thm}
\label{thm:quasihopfish} Let $(A,\Delta,\epsilon, \Phi)$ be a
biunital quasi-bialgebra.  Then  $(A\otimes  A)/W$, where $W$ is
the left ideal generated by
 $$\{\epsilon(a)(1\otimes 1)-\Delta(a)|a\in A\},$$ is a preantipode for the modulation of $A$.
  If $A$ is a quasi-Hopf algebra, with antipode $(S, \alpha, \beta)$, then $(A\otimes  A)/W$
  is isomorphic to the modulation  $\bS$, and $(A,\bDelta,\bepsilon,\bS)$
 is a hopfish algebra.
\end{thm}
\section{Morita invariance}
\label{sec-morita}

The following theorem shows that, with our definition of hopfish
algebra, we are on the right track toward defining a Morita
invariant notion.

\begin{thm}
\label{thm:morita-hopf} Let $A$ be a quasi-Hopf algebra and $B$
an algebra Morita equivalent to $A$. Then $B$ is a sesquiunital
sesquialgebra with a preantipode.
\end{thm}
$\pf$ Let $P$ be an $(A,B)$-bimodule, and $Q$ a $(B,A)$-bimodule,
inverse to one another in the category \Alg.  We recall the
hopfish structure on $A$ defined in Theorem
\ref{thm:quasihopfish}, with
\[
\bepsilon^A=k,\ \ \ \bDelta^A=A\otimes A,\ \ \ \bS^A=A\otimes A/W.
\]

We use the bimodules $P$ and $Q$ to define:
\[
\begin{array}{c}
\bepsilon^B:=\bepsilon^A\otimes_A P,\ \ \ \ \
\bDelta^B:=\big(Q\otimes Q\big)\otimes_{A\otimes A}\bDelta^A
\otimes_A P,
\end{array}
\]
It is straightforward to check that these form a sesquiunital
sesquilinear algebra structure on $B$.

Now we define
\[
\bS^B:=\big( Q\otimes Q\big)\otimes_{A\otimes A}\bS^A.
\]
\begin{rmk}
{\rm We remark that our definition of the antipode $\bS^B$ only uses the bimodule $Q$, not P. This is because $Q$ is a $(B, A)$ bimodule, and therefore is also an $(A^{op}, B^{op})$ bimodule naturally. Since $\bS^A$ is an $(A, A^{op})$ bimodule, $Q\otimes _A \bS^B\otimes_{A^{op}} Q$ defines a $(B, B^{op})$ bimodule, which is isomorphic to $(Q\otimes Q)\otimes _{A\otimes A}\bS^A$.}
\end{rmk}

In the following, we will show that $\bS^B$ is a preantipode, i.e.
\[
\Hom_k(k, \bS^B )\cong \Hom_B(\bepsilon^B, \bDelta^B).
\]

According to our definitions, we have
\[
\Hom_B(\bepsilon^B, \bDelta^B)=\Hom_B\big(\bepsilon^A\otimes_A
P,\big( Q\otimes Q\big)\otimes_{A\otimes A}\bDelta^A
\otimes_AP\big).
\]

Since the Morita equivalence between $A$ and $B$ defines an
equivalence of right-module categories, we have a natural
isomorphism
\[
\Hom_B\big(\bepsilon^A\otimes_A P,\big( Q\otimes
Q\big)\otimes_{A\otimes A}\bDelta^A\otimes_AP\big)\cong
\Hom_A\big(\bepsilon^A, \big(Q\otimes Q\big)\otimes_{A\otimes
A}\bDelta^A \big).
\]

The space  $\Hom_A\big(\bepsilon^A, \big(Q\otimes
Q\big)\otimes_{A\otimes A}\bDelta^A\big)$ consists of $k-$linear
morphisms from $\big(Q\otimes Q\big)\otimes_{A\otimes A}\bDelta^A$
to $k$, vanishing on the $A-$submodule $\widetilde{W}$ spanned by
\[
(q_1\otimes q_2)\otimes_{A\otimes A}(a_1\otimes
a_2)(\epsilon(a)1\otimes 1-\Delta(a)),\ q_i\in Q,\ a,\ a_i\in A,\
i=1,2.
\]

The $A-$submodule $\widetilde{W}$ is isomorphic to $(Q\otimes
Q)\otimes_{A\otimes A}W$, where $W$ is defined as in Theorem
\ref{thm:quasihopfish}. Therefore, $\Hom_A\big(\bepsilon^A,
\big(Q\otimes Q\big)\otimes_{A\otimes A}\bDelta^A\big)$ is
isomorphic to the $k$-dual of the quotient
\begin{equation}
\label{eq:dual-anti} (Q\otimes Q)\otimes_{A\otimes
A}\bDelta^A/\widetilde{W}\cong (Q\otimes Q)\otimes_{A\otimes
A}(A\otimes A/W).
\end{equation}
Replacing $A\otimes A/W$ by $\bS^A$ in (\ref{eq:dual-anti}), we
have
\[
\Hom_B(\bepsilon^B, \bDelta^B)\cong \big((Q\otimes
Q)\otimes_{A\otimes A}\bS^A\big)^*\cong (\bS^B)^*.
\]

Therefore, $\bS^B$ defines a preantipode on $(B, \bDelta^B,
\bepsilon^B)$. $\Box$

\bigskip
Now we study when the  sesquiunital sesquialgebra just defined is
a hopfish algebra, i.e. when  $\bS^B$ is isomorphic to $B$ as a
left
 $B$-module.

We introduce the following special type of module over a hopfish
algebra.
\begin{dfn}\label{dfn:self-conjugate}
Let be $A$ be a hopfish algebra with  antipode bimodule $\bS$,
and let $X$ be a right $A$-module and therefore a left
$A^{op}$-module. Then $X$ is {\bf self-conjugate} if $\Hom_A(A,
X)$ is isomorphic to $\bS \otimes_{A^{op}} X$ as a left
$A$-module.
\end{dfn}
\begin{rmk}
{\em We remark that the category of finite dimensional left modules over
a quasi-Hopf algebra is a rigid monoidal category. A self-dual
module $X$ of a quasi-Hopf algebra $A$ is a
self-dual object in the category of finite
dimensional modules, i.e. $\Hom_k(k, X)$ is isomorphic to
$\bS\otimes_{A^{op}} X$. }
\end{rmk}

We can understand the definition of a self-conjugate module
geometrically as follows. A hopfish algebra $A$ can be thought as
functions on a ``noncommutative space with group structure" $G$.
If we view a finite projective right $A$-module $X$  as the space
of sections of a ``vector bundle" $E$ over $G$,  $\Hom_A(A, X)$
corresponds to the space of sections of the dual bundle $E^*$, and
$\bS\otimes_{A^{op}}X$ is the pullback of the bundle $E$ by the
``inversion" map $\iota$ of $G$. The self-conjugacy condition on
$E$ says that $E^*$ is isomorphic to $\iota^*{E}$.

\begin{prop}
\label{prop-morita}With the same assumptions and notation as in
Theorem \ref{thm:morita-hopf}, if the $(B,A)$-Morita equivalence
bimodule $Q$ is self-conjugate as a right $A$-module, then $B$ is
a hopfish algebra with antipode $\bS^B$ defined in Theorem
\ref{thm:morita-hopf}.
\end{prop}
\pf Recall that the preantipode on $B$ defined in Theorem
\ref{thm:morita-hopf} is equal to
\[
(Q\otimes Q)\otimes _{A\otimes  A}\bS^A.
\]
Since $Q$ is a right $A-$module, it is also a left
$A^{op}-$module, and the preantipode $\bS ^B $ can be rewritten as
\[
Q\otimes_A \bS^A\otimes_{A^{op}} Q.
\]

Since $Q$ is self conjugate, we have
\[
\bS^A\otimes_{A^{op}} Q\cong \Hom_A(A, Q),
\]
and so
\[
Q\otimes_A \bS^A\otimes_{A^{op}} Q\cong Q\otimes_A  \Hom_A(A,Q).
\]

When $Q$ is a Morita equivalence bimodule between $A$ and $B$, $Q$
is a finitely generated projective $A$-module and $B\cong
\Hom_A(Q, Q)=Q\otimes_A\Hom_A(A, Q)$. This shows that $Q\otimes_A
\bS^A\otimes_{A^{op}} Q$ is isomorphic to $B$ as a left
$B$-module. $\Box$

The following example is a special case of Proposition
\ref{prop-morita}. We remark that given a (quasi)-Hopf algebra
$A$, the matrix algebra $M_n(A)$ of $n\times n$ matrices with
coefficients in $A$ is not a (quasi-)Hopf algebra when $n\geq 2$.
\begin{ex}
\label{ex-matrices} Let $A$ be a quasi-Hopf algebra with
$\bepsilon^A=k$, $\bDelta^A=A\otimes A$, and $\bS^A=A$. Then the
$n\times n$ matrix algebra $M_n(A)=B$ with coefficients in $A$ is
a hopfish algebra. We consider $Q=A^n$ as a space of column
vectors, so that it has the structure of an $(M_n(A),A)$-bimodule,
The counit $\bepsilon^B$ is $A^n$ viewed as row vectors, i.e. as a
$(k, M_n(A))$-bimodule. The coproduct $\bDelta^B$ is isomorphic to
\[
\big( A^{ n}\otimes A^{n} \big)\otimes_{A\otimes
A}(A\otimes A)\otimes_A (A^{n})^T=\big( A^{n}\otimes
A^{n} \big)\otimes_{A\otimes A} (A^{n})^T.
\]
$\bS^B$  is equal to  $A^{ n}\otimes_A A\otimes_{A^{op}}A^{n}$.
$A^{n}\otimes_A A \otimes_{A^{op}} A^{n}$ is isomorphic to
$M_n(A)$ as an $(M_n(A),{M_n(A)}^{op})$-bimodule, where the left
$M_n(A)$ module structure is from the standard left
multiplication, while the right $M_n(A)^{op}$ module structure is
the composition of the left multiplication, transposition of
matrices, and the antipode on $A$. Therefore, $B=M_n(A)$ is a hopfish
algebra.
\end{ex}

The following example shows that the self-conjugacy condition in
Proposition \ref{prop-morita} can not be eliminated.
\begin{ex}
{\em We consider the cyclic group $\integers/3\integers$ with
elements $0, 1, 2$. The algebra $A$ of functions on
$\integers/3\integers$ is a commutative Hopf algebra spanned by
the characteristic functions
 $e_0$, $e_1$, and $e_2$.
We notice that the  $e_i$'s are projections in $A$, and denote the submodule
$e_i A$ by $A_i$.   Now consider
the following projective module over $A$
\[
Q=A_0^r\oplus A_1^s\oplus A_2^t,
\]
where $r,s,t$ are nonnegative integers.  Then
$$B=\Hom_A(Q,Q)=A_0^{r^2}\oplus A_1^{s^2}\oplus
A_2^{t^2}.$$
 It is not difficult to see that $Q$ is self-conjugate if and only if $s=t$.

We calculate the expression for $\bS^B$ in Theorem
\ref{thm:morita-hopf} as follows,
\[
\begin{split}
&\; (Q\otimes Q)\otimes_{A\otimes A}\bS^A\\
 &=\big(A_0^{r}\oplus A_1^{s}\oplus A_2^{t})\otimes
 (A_0^{r}\oplus A_1^{s}\oplus
A_2^{t})\big)\otimes_{A\otimes A} \bS^A\\
&=\big(A_0^{r}\otimes (A_0^{r}\oplus A_1^{s}\oplus
A_2^{ t})\big)\otimes _{A\otimes
A}\bS^A\\
&\oplus \big(A_1^{s}\otimes (A_0^{r}\oplus
A_1^{s}\oplus A_2^{t})\big)\otimes _{A\otimes
A}\bS^A\\
&\oplus \big(A_2^{t}\otimes (A_0^{r}\oplus
A_1^{s}\oplus A_2^{t})\big)\otimes _{A\otimes A}\bS^A.
\end{split}
\]

We look at the tensor product $(A_i\otimes A_j)\otimes_{A\otimes
A}\bS^A$. By Theorem \ref{thm-hopfhopfish}, the antipode bimodule
$\bS^A$ is isomorphic to $A$. Therefore $(A_i\otimes
A_j)\otimes_{A\otimes A}\bS^A$ is equal to
\[
(A_i\otimes A_j)\otimes_{A\otimes A}A=A_i\otimes_A A_j,
\]
where the left $A$-module structure on $A_j$ is the composition
of the right multiplication with the antipode map $S:A\to A$.

We notice that $S(e_i)e_j=0$ if $S(e_i)\ne e_j$. Therefore,
\[
A_i\otimes_A A_j=\left\{\begin{array}{ll}0&S(e_i)\ne e_j,\\
A_i&S(e_i)=e_j\end{array} \right .
\]

We conclude that $\bS^B=A_0^{ r^2}\oplus A_1^{st}\oplus
A_2^{st}$. We observe that $\bS^B$ is isomorphic to $B$ as a
left $B$ module if and only if $s=t$.

Therefore, $\bS^B$ is isomorphic to $B$ if and only if $Q$ is a
self-conjugate $A$-module. }
\end{ex}

We define a notion of  Morita equivalence between hopfish
algebras.
\begin{dfn}
\label{dfn:hopfish-morita}Let $(A, \bepsilon^A, \bDelta^A, \bS^A)$
and $(B, \bepsilon^B, \bDelta^B, \bS^B)$ be two hopfish algebras.
Then $(A, \bepsilon^A, \bDelta^A, \bS^A)$ is Morita equivalent to
$(B, \bepsilon^B, \bDelta^B, \bS^B)$ if there there is an
$(A,B)$-bimodule $_AP_B$ and a $(B,A)$-bimodule $_BQ_A$ satisfying
\begin{enumerate}
\item $P\otimes _B Q=A$, and $Q\otimes _A P=B$.
\item $\bepsilon^B=\bepsilon^A\otimes_A P$,
\item $\bDelta^B=(Q\otimes Q)\otimes_{A\otimes A}\bDelta^A \otimes_A
P$,
\item $\bS^B=(Q\otimes Q)\otimes_{A\otimes A}\bS^A$.
\end{enumerate}
\end{dfn}

\begin{prop}
Definition \ref{dfn:hopfish-morita} defines an equivalence
relation among hopfish algebras.
\end{prop}
\pf Straightforward check. $\Box$
\section{Hopfish structures on $k^n$}\label{sec:kn}
In this section, we give examples of hopfish algebras which are
not Morita equivalent to  modulations of Hopf algebras.
In particular, we will describe hopfish structures on
the commutative algebra $k^G$ of $k$-valued functions on a finite set $G$ which do not correspond to group structures on $G$.

We
may identify the $r$-th tensor power of $k^G$ with $k^{G^r}$.
Since this algebra is commutative, we can also identify
$(k^g)^{op}$ with $k^g$.

If $G$ is a semigroup, $k^G$ is a bialgebra with coproduct
$\Delta(a)(g,h) =a(gh)$, with a counit $\epsilon(a)=a(e)$ when $G$
has an identity element $e$.  When $G$ is a group, we get a Hopf
algebra structure by setting $S(a)(g)=a(g^{-1})$.

Now let $G$ be a groupoid.  We may make the same definitions as
above, adding that $\Delta(a)(g,h) $ should be $0$ when $gh$ is
not defined, and $\epsilon(a)$ is the sum of the values of $a$ on
all the identity elements. When $G$ is not a group, $k^G$ is
no longer a Hopf algebra, but rather a {\bf weak Hopf algebra}
(Example 2.3, \cite{n:weak-hopf}), since $\Delta$ is not unital and
$\epsilon$ is not even an algebra homomorphism. When $G$ is a
groupoid, we have two algebra morphisms $\alpha$,
$\beta$:$k^{G^0}\to k^G$ as the lifts of the source and target
maps. The coproduct $\Delta$ is defined on
$k^G\otimes_{k^{G^0}}k^G$ by $\Delta(a)(g,h)=a(gh)$, and counit
$\epsilon:k^G\to k^{G^0}$ by $\epsilon(a)(e)=a(e)$, and the
antipode $S$ is defined by $S(a)(g)=a(g^{-1})$. $(k^G, \alpha,
\beta, \Delta, \epsilon, S)$ is a {\bf quantum groupoid},
\cite{lu:quant}. It turns out that we can still form the
modulation of these operators, and we still get a hopfish algebra
because of the commutativity of the algebras $k^G$ and $k^{G^0}$. To
prove this, we will look at a more general situation.

Any sesquialgebra coproduct on $A=k^G$ is an $(A\otimes
A,A)$-bimodule, i.e. a module $\bDelta$ over $k^{G\times G\times
G}$.  Such a module decomposes into submodules supported at the
points of $G^3$.  For our purposes, we will assume that these are
free modules of finite rank.  Then $\bDelta$ is determined up to
isomorphism by the dimensions $d^{g}_{hk}$ of the components
$\bDelta^g_{hk}$, for $(g,h,k)\in G^3$.  It is straightforward to
check that the condition for coassociativity is precisely that the
$d^{g}_{hk}$'s be the structure constants of an associative
algebra $A'=\integers^G$ over $\integers$, i.e. a ring. Namely, identifying each element of $G$ with its characteristic function, we have
$gh = \sum_k d_{gh}^k k$.  Similarly, a $(k,A)$-bimodule $\bepsilon$ with free submodules $\bepsilon^g$ as
components is determined by the dimensions $e^g$ of $\bepsilon^g$,
and this bimodule is a counit precisely when $e:= \sum_g e^g g$
is a unit element for $A'$. We call such sesquiunital
sesquialgebras of {\bf finite free type}. Thus we have shown:

\begin{prop}
There is a one to one correspondence between sesquiunital sesquialgebra structures of finite free type on
$k^G$ and unital ring
structures on $\integers^G$  for which the structure constants and
the components of the unit are nonnegative.
\end{prop}

The best known examples of such rings  are the
monoid algebras.  If $G$ is a monoid, then we may define $\delta_{hk}^g$
to be the characteristic function of the graph $g=hk$ of
multiplication and $e^g$ to be the characteristic function of
the identity element. The corresponding sesquialgebra is
just the modulation of the dual to the monoid bialgebra $A'$.

With this example in mind, we may think of a general structure of
convolution type on $\integers^G$ as corresponding to a
``product'' operation on $G$ in which the product of any two
elements is a (possibly empty) subset of $G$ whose elements are
provided with positive integer ``multiplicities'.  We will call
such a subset a ``multiple element''; the identity is also such a
multiple element.  (Of course, any ring structure may be viewed
in this way, if we allow the multiplicities to be arbitrary
integers).

To begin our analysis of these structures, we show that there are
restricted possibilities for the unit.

\begin{prop}
\label{prop-counit} The $e^g$'s are either $0$ or $1$.
\end{prop}

\pf  Given $g$, by the counit property,
\[
\sum_{k}e^k d^{g}_{gk}=\delta_{gg}=1,
\]
we have that at there is at least one $k\in G$, such that
$d^g_{gk}\ne0$.

 By the counit property again, we
have
\[
e^g \leq e^g d^g_{gk} \leq \sum_{h}e^h d^g_{hk}=\delta_{gk}\leq 1.
\]
 \qed

We will denote by $G^{0}$ the support of the unit, i.e. the set of
$g\in G$ for which $e^g=1$.  This set will play the role of
identity elements in
$G$.

As long as $G$ is nonempty, so is $G^{0}$.  In fact, we have the
following:

\begin{prop}\label{prop-unit}
Given any $g$ in $G$, there are unique elements $l(g)$ and $r(g)$
in $G^{0}$ such that, for all $h\in G^{0}$, $d_{hg}^k =
\delta_{hl(g)} \delta_{gk}$ and $d_{gh}^k =
\delta_{r(g)h}\delta_{gk}$.
\end{prop}

\pf This is again a straightforward corollary of the counit
property.
\[
\begin{array}{ll}
\sum_{k}e^k d_{gk}^h=\delta_{gh},& \sum_{g}e^g
d_{gk}^h=\delta_{kh}.
\end{array}
\]
Therefore $\sum_{g\in G^0} d_{gh}^k=\delta_{kh}$. So $d_{gh}^k=0$
when $k \neq h$ and there exists a unique element  $g_0\in G^0$
such that $d_{gh}^h$ equals 1 for $g=g_0$ and 0 for all other $g$. We let $l(h)$ be
this $g_0$. So the first equation is proven, and the second is proven
by a similar argument. \qed

Since the sum of the elements of $G^{0}$ is the unit of $k^G$, it
is idempotent, from which it follows that $k^{G_{0}}$ is a
subalgebra.  In fact, one may show:

\begin{prop}
The elements of $G^{0}$ form a set of orthogonal idempotents in
$\integers^G$.  In other words, the algebra structure on the subalgebra
$\integers^{G^{0}}$ of $A'$ is just pointwise multiplication.
\end{prop}

\pf This follows from uniqueness in Proposition \ref{prop-unit}.
\qed

We also have:

\begin{prop}\label{prop-lr}
For all $g$ and $h$ in $G$, if $d^k_{gh}\ne0$, $l(k)=l(g)$ and
$r(k)=r(h)$. If $r(g)$ is not equal to $l(h)$, then $gh=0$ in $G$.
In particular, $l(h)=h=r(h)$ for all $h\in G^0$.
\end{prop}

\pf Coassociativity gives us
\[
\sum_{s}d^k_{l(g)s}d^s_{gh}=\sum_{s}d^s_{l(g)g}d^k_{sh}.
\]

By Proposition \ref{prop-unit}, $d^s_{l(g)g}=\delta_{gs}$.
Therefore, the right hand side of the equation is equal to
$d^k_{gh}\ne0$.

On the left hand side, according to Proposition \ref{prop-unit},
 $d^k_{l(g)s} \neq 0$ only  when
$l(s)=l(g)$ and $k=s$. Therefore, if $d^k_{gh}\ne0$,then
$d^k_{l(g)k}=1$, so $l(k)=l(g)$.
Similar arguments  show that $r(k)=r(h)$.

If $r(g)\ne l(h)$, again by coassociativity, we have
\[
d^k_{gh}=\sum_{s}d^s_{gr(g)}d_{sh}^k=\sum_{s}d^k_{gs}d^s_{r(g)h}.
\]
According to Proposition \ref{prop-unit}, if $r(g)\ne l(h)$,
$d^s_{r(g)h}=0, \forall s$; therefore, $d^k_{gh}=0$. \qed

We now have retractions $l$ and $r$ from $G$ onto $G^{0}$ which are
like the ``target'' and ``source'' maps from a category to its set
of identity elements.   In fact,
in terms of the multiplicative structure on $G$ corresponding to
the algebra structure on $A'$, we have $l(g)g=gr(g)=g$; in
particular, these products are single valued and without
multiplicities.  We might call
 $G$ a ``hypercategory''.   The
composition of morphisms is a ``multiple morphism'' between two
definite objects.

We will show next that, when $k^G$ has an antipode and is hence a
hopfish algebra, the underlying multiplicative structure on $G$ has
inverses and the property that $gh$ is nonzero whenever
$r(g)=l(h)$.  We will call such a structure a
``hypergroupoid'' (see Definition \ref{dfn-hyper-gpd}).\footnote{The notion of group with multiple-valued multiplication has a long history.  We give just one reference, to the 1939 paper by Kuntzmann \cite{ku:contribution}, which contains references to even earlier work.}

According to Definition \ref{dfn-antipode}, an antipode is a
$(k^G,k^G)$-bimodule $\bS$ whose dual is isomorphic to
$\Hom_{k^G}(\bepsilon, \bDelta)$.

We recall the definition of $\bepsilon$ and $\bDelta$
\[
\bepsilon=\oplus_g \bepsilon^g,\ \ \ \
\bDelta=\oplus_{g,h,t}\bDelta^t_{gh}.
\]

Therefore, $\Hom_{k^G}(\bepsilon, \bDelta)$ may be written as
\[
\big(\oplus_s \bepsilon^s\big)\otimes
_{k^G}\big(\oplus_{g,h,t}\bDelta^t_{gh}\big)^*.
\]

We notice that $k^G$ acts via the upper indices of $\bepsilon$ and
$\bDelta$ by componentwise multiplication. Therefore, the above
expression for $\Hom_{k^G}(\bepsilon, \bDelta)$ can be simplified
to
\[
\oplus_{g,h}\big(
\oplus_{t}{\bepsilon^t}^*\otimes\bDelta^t_{gh}\big)^*,
\]
which is isomorphic to
\[
\oplus_{g,h} \big(\oplus_t
\Hom_k(k, {\bepsilon^t}^*\otimes\bDelta^t_{gh})\big)\cong
\Hom_k\big(k, \oplus_{g, h}(\oplus_t
{\bepsilon^t}^*\otimes\bDelta^t_{gh})\big).
\]

Therefore, $\bS$ is isomorphic to $\oplus_{g,h}(\oplus_t
{\bepsilon^t}^*\otimes\bDelta^t_{gh})$ as a $(k^G,k^G)$ bimodule.

When $\bS$ is an antipode,  $\bS$ is by definition isomorphic to $k^G$ as a
left $k^G$-module. Therefore, if we write $\bS$ as $\oplus_{g,h}\bS_{gh}$, for any
fixed $g$ there exists a unique element $h\in G$ such that
$\dim(\bS_{gh})=1$, and $\dim(\bS_{gh'})=0$ for  $h' \neq h$. Hence,
we may define a map $\sigma: G\to G$ by $g\mapsto h$.

Since  $\Hom_k(k, \bS)\cong\Hom_{k^G}(\bepsilon, \bDelta)$, we know
that
\[
\dim(\oplus_t {\bepsilon^t}^*\otimes\bDelta^t_{gh})=\delta_{\sigma(g)h}, i.e.
\]
\begin{equation}
\label{eq-antipode} \sum_t e^t d^t_{gh}= \delta_{\sigma(g)h}.
\end{equation}

From this, we have the following:
\begin{prop}
For any $g\in G$, define $g^{-1}$ to be $\sigma(g)$. Then
there is a unique $s\in G^0$ such that
\[
d_{gg^{-1}}^s=1.
\]
In fact,   $s=l(g)=r(g^{-1})$.
\end{prop}
\pf By Equation (\ref{eq-antipode}), we have
\[
\sum_t e^t d^t_{gg^{-1}}= \delta_{\sigma(g)g^{-1}}=1.
\]
Therefore, there is a unique element $s\in G^0$ such that
$d^{s}_{gg^{-1}}=1$, and $d^t_{gg^{-1}}=0$ for all other $t\in
G^0$. By Proposition \ref{prop-lr}, $s=l(g)=r(g^{-1})$. \qed

We also have another characteristic property of categories (though here we can only prove it in the presence of an antipode).
\begin{prop}\label{prop-product}
If $r(g)=l(h)$, then there exists $s\in G$, such that $d_{gh}^s\ne0$.
\end{prop}
\pf Using coassociativity, we have
\[
\sum_s d^g_{gs}d^s_{hh^{-1}}=\sum_{s}d^s_{gh}d^g_{sh^{-1}}.
\]
Since $d^{l(h)}_{hh^{-1}}=1$ and $d^g_{gl(h)}=1$(since
$r(g)=l(h)$), the left hand side of the above equation is not
equal to 0. This implies that, on the right hand side, there is at
least one term which is not equal to 0. Therefore, there exists $s\in
G$, such that $d^s_{gh}\ne 0$. \qed

\begin{question}
Is the inversion map $\sigma:G\to G$  involutive?
\end{question}

To summarize  the arguments above, we define
the ``combinatorial" objects  associated to  hopfish
structures on $k^G$:
\begin{dfn}
\label{dfn-hyper-gpd}A  {\bf hypergroupoid} is  a set $G$ with the
following operations $(\cdot, l,r, ^{-1})$.
\begin{enumerate}
\item There is a multi-valued associative binary operation $\cdot$ on
$G$. More precisely, $\forall g, h\in G$, $\ g\cdot h$ is an
element
  of $\integers_+^G$, where $\integers_+$ is the semiring of
  nonnegative integers.   When
this product is linearly extended to a product on $\integers_+
^G$, we have $g\cdot (h\cdot k)=(g\cdot h)\cdot k$.
\item There is a subset $G^0\subset G$ with maps
$l,r:G\to G^0$ such that $l(g)\cdot g=g\cdot r(g)=g$, $\forall
g\in G$. The product of $g$ and $h$ is nonzero if and only if
$r(g)=l(h)$.
\item There is an inverse operation $g\mapsto g^{-1}$ on $G$ such that $g\cdot g^{-1}|_{G^0}= l(g)$ $\forall g\in
G$, and if $h\ne g^{-1}$, then $g\cdot h |_{G^0}=0$.
\end{enumerate}
\end{dfn}

Note that the inverse operation is
determined by the product operation and $G^0$.

We now look at the commutative algebra $k^G$ of $k$-valued
functions on a hypergroupoid $G$. The duals of the maps $l,
r:G\to G^0$ define morphisms from $k^G$ to $k^{G^0}$. The multiplication on $G$ defines a (nonunital)
coproduct homomorphism $k^{G\times G}\leftarrow k^G$ whose modulation is a coproduct bimodule, the embedding map from $G^0$ to $G$ makes
$k^{G^0}$ into a counit bimodule, and the inversion map defines an
antipode. These define a hopfish algebra structure on $k^G$.

\begin{thm}
The  hopfish structures of free finite type on $k^G$ are in one to
one correspondence with the hypergroupoid structures on $G$.
\end{thm}

\begin{ex}
\label{ex:noar} {\em Here is an example of a hypergroupoid which
is not a groupoid, based on  Example 8.19 in  \cite{ent:fusion}.\footnote{The hypergroupoid itself , when $n=1$, already appears as the first example in \cite{ku:contribution}!}
 Let
$G=\{ e, g\}$, with multiplication and inversion given by
\[ eg=ge=g, \quad ee=e, \quad gg=e+ng, \quad e^{-1}=e, \quad g^{-1}=g \]
where $n$ is a nonnegative integer. $G^0=\{e\}$ and $l(g)=r(g)=e$. The
algebra $A'$ associated to this hypergroupoid  is $\integers
[x]/\{x^2=1+nx\}$. The corresponding hopfish algebra $k^G$ is not
a quasi-Hopf algebra when $n=1$ and $k$ a field. We explain the
reason  below.

The hopfish algebra structure of $k^G$ is in fact a weak Hopf
algebra, with $\epsilon(\ e)=1$, $\epsilon(\ g)=0$, $\Delta(\ e)=\
g\otimes\ g+\ e\otimes\ e$ and $\Delta(\ g)=\ e\otimes\ g+\
g\otimes\ e+\ g\otimes\ g$. Since a $k^G$ module can be decomposed
into submodules supported at points of $G$, the representation
ring of $k^G$ is generated by two elements 1 and $X$ corresponding
to the 1-dimensional $k^G$ module supported at $e$ and $g$
respectively. Since $k$ is a field, 1 and $X$ are just
1-dimensional $k$-vector spaces. Using the formulas for the
coproduct and counit, it is easy to check that this representation
ring is the Grothendieck ring of what is called Yang-Lee category
in \cite{os:fusion}, namely it is generated by 1 and $X$ with the
relation $X\otimes X=1\oplus X$. The Frobenius-Perron dimension of
the element $1$ is 1, while the Frobenius-Perron dimension of the
element $X$ is the irrational number $(1+\sqrt{5})/2$. According
to Theorem 8.33, \cite{ent:fusion}, the Frobenius-Perron dimension
of any finite dimensional module over a finite dimensional
quasi-Hopf algebra must be a positive integer, which is equal to
the dimension of the module.  This shows that $k^G$ is not Morita
equivalent to a quasi-Hopf algebra. }
\end{ex}

\def\cprime{$'$}

\noindent{Xiang Tang, xtang@math.ucdavis.edu}\\
\noindent{\rm  Department of Mathematics, University of
California, Davis, CA, 95616, USA. } \vskip 3mm

\noindent{Alan Weinstein, alanw@math.berkeley.edu}\\
\noindent{\rm Department of Mathematics, University of California,
Berkeley, CA, 94720, USA.}\vskip 3mm

\noindent{Chenchang Zhu, zhu@math.ethz.ch}\\
\noindent{\rm Department Mathematik, Eidgen\"ossische Technische
Hochschule, 8092 Z\"urich, Switzerland.}

\end{document}